\documentclass[12pt]{article}
\usepackage{amssymb}
\usepackage{amsmath}
\usepackage{relsize}
\usepackage{theorem}
\usepackage{latexsym}
\newtheorem{thm}{Theorem}[section]
\newtheorem{prop}{Proposition}[section]
\newtheorem{lem}{Lemma}[section]
\newtheorem{cor}{Corollary}[section]
\newtheorem{de}{Definition}[section]
\newtheorem{ex}{Example}[section]
\title{Convergence of Infinite Composition of Entire Functions}
  \author{Shota Kojima}
  \date{}
\begin{document}
\def\R{\mathop{\mathcal{R}}}
\def\L{\mathop{\mathcal{L}}}
\maketitle
\begin{abstract}
The purpose of the present article is to obtain the condition that the function defined by infinite composition of entire functions  becomes an entire function. Moreover, as an example of such functions, we study a function called Poincar\'e function.
\end{abstract}
\section{Introduction}
It seems that there are no article studying infinite composition of functions with a similar purpose to this article.
In order to state our theorem, we require the following notation.
\begin{de} Let $f\circ g$ be the composition of functions $f$ and $g$, that is,
\[ (f \circ g)(z) := f(g(z)) . \]
We denote $(f \circ g)(z)$ by $f(z) \circ g(z)$ for convenience of expression.
\end{de}
For example 
\[
(z+1) \circ (z+2) = z+3.
\]
\begin{de} Given integers $d,N$ with $N\geq d$, define
\begin{eqnarray*}
\R_{n=d}^{N}f_n(z) : &=& f_d(z) \circ f_{d+1}(z) \circ f_{d+2}(z)  \circ \cdots \circ f_N(z) \\
                     &=& f_d(f_{d+1}( \cdots f_{N-1}(f_N(z)) \cdots).
\end{eqnarray*}
\end{de}
Some important functions, such as $\sin z , e^z$, are expressed by infinite composition of polynomials as follows:
\begin{prop}\label{specialvalue}
For any $z \in \mathbb{C}$, we have
\begin{eqnarray*}
\frac{1}{2}(e^{2z}-1) &=& \R_{n=1}^{\infty} \left( z+\frac{z^2}{2^n} \right)    \, ,\\
\sin \left( \frac{2z}{\sqrt{3}}+\frac{\pi}{6} \right)-\frac{1}{2} &=& \R_{n=1}^{\infty} \left( z+\frac{z^2}{(-2)^n} \right) \, ,\\
\left( \sinh \sqrt{z} \right)^2 &=& \R_{n=1}^{\infty} \left( z+ \frac{z^2}{4^n} \right) \, .
\end{eqnarray*}
\end{prop}
These equalities are proved in Section $3$.
From the equalities above, we expect that there are remarkable functions defined by infinite composition of entire functions.
Thus it is significant to study the convergence of
\[
\R_{n=1}^{\infty}f_n(z)=\lim_{N \rightarrow \infty} \R_{n=1}^{N} f_n(z) = f_1(z) \circ f_2(z) \circ \cdots , \,\, \text{where} \,\,f_n(z) \,\, \text{is entire}.
\]
Our main purpose is to prove the following theorem.
\begin{thm}\label{entire1} 
Let $c_{n,r} \, (n=1,2,\ldots, r=2,3, \ldots) $ be complex numbers such that 
\[
f_n(z) := z+ \sum_{r=2}^{\infty} c_{n,r} z^r 
\]
are entire functions. 
We set 
\[
 C_n := \max_{r=2,3,4,\ldots} \{ |c_{n,r}|^{1/(r-1)} \}. 
\]
Suppose that the series 
\[
\sum_{n=1}^{\infty} {C_n}
\]
is convergent. Then the sequence of functions
\[
 \R_{n=1}^{N} f_n(z) = \R_{n=1}^{N} \left( z + c_{n,2} z^2 + c_{n,3}z^3 + \cdots + c_{n,p} z^p + \cdots \right) 
\]
is uniformly convergent on arbitrary closed disk. 
In particular, the limit function 
\begin{eqnarray*}
\R_{n=1}^{\infty} f_n(z) &=& \R_{n=1}^{\infty} \left( z + c_{n,2} z^2 + c_{n,3}z^3 + \cdots + c_{n,p} z^p +\cdots \right) \\
                         &= &\lim_{N \rightarrow \infty}\R_{n=1}^{N} \left( z + c_{n,2} z^2 + c_{n,3}z^3 + \cdots + c_{n,p} z^p + \cdots \right) 
\end{eqnarray*}
is entire.
\end{thm}
Considering the case where $f_n(z)$ is a polynomial of degree $2$, we obtain the following.
\begin{cor}
Let $\{c_n\}_{n=1}^{\infty}$ be a sequence of complex numbers such that
\[
\sum_{n=1}^{\infty} |c_n| 
\]
is convergent. Then the function
\[
\R_{n=1}^{\infty} \left( z + c_nz^2 \right) = \lim_{N \rightarrow \infty} \R_{n=1}^{N} \left( z + c_nz^2 \right) 
\]
is uniformly convergent on every compact subset of $\mathbb{C}$, and it defines an entire function.
\end{cor}

\begin{ex}{\rm
The infinite composition
\[
\R_{n=1}^{\infty} \left( z + \frac{z^3}{n^3} \right)
\]
is an entire function. Indeed, for $ c_{n,3} = n^{-3}$, the series $\sum_{n=1}^{\infty}|c_{n,3}|^{1/2}= \sum_{n=1}^{\infty} n^{-3/2}$ is convergent.}
\end{ex}
\begin{ex}{\rm{
Let $s$ be a complex number with $|s|>1$. 
Then the infinite composition
\[
F(z) :=\R_{n=1}^{\infty} \left( z+ \frac{z^2}{s^n} \right) 
\]
is an entire function. Indeed, for $ c_{n,2} = s^{-n}$, the series $\sum_{n=1}^{\infty}|c_{n,2}|=\sum_{n=1}^{\infty}|s|^{-n} $
is convergent. This function is studied in Section 3.}}
\end{ex}
We now introduce Poincar\'e functions. (For more details, see \cite{Gregory}.) The meromorphic functions $f(z)$ satisfying the following functional equation are called Poincar\'e functions (\cite{Valiron}):
\[
f(sz) =h(f(z)),
\]
where $s$ is a complex number with $|s|>1$, and $h(z)$ is a rational function.
The function $F(z)$ in Example $2$ satisfies
\[
F(sz) = s F(z) + s F(z)^2
\]
(see Section $3$).
Thus the function $F(z)$ can be regarded as a Poincar\'e function. 
Poincar\'e functions have been studied by some mathematicians (\cite{Gregory}, \cite{Mustafa}, \cite{Poincare}, \cite{Valiron}).
However it seems that the expression of Poincar\'e functions by $\R$ is not known.

\section{Proof of Theorem \ref{entire1} }
In this section we shall give a proof of Theorem \ref{entire1}.
First we define
\begin{de}
For any analytic function
\[
f(z) = \sum_{n=0}^{\infty} a_n z^n,
\]
we define 
\[
\hat{f}(z) =  \sum_{n=0}^{\infty} |a_n| z^n.
\]
\end{de}
Second, we prove lemmas needed later. 
\begin{lem}\label{ineq:firstlem1}
Let 
\begin{eqnarray*}
f(z) : = z+ \sum_{n=2}^{\infty} a_n z^n, \quad  g(z) : = z+ \sum_{n=2}^{\infty} b_n z^n
\end{eqnarray*}
be entire functions.
Then for every $z \in \mathbb{C}$,
\begin{eqnarray}
{\widehat{ f \circ g }} (|z|) &\leq & \hat { f} (\hat{g}(|z|) ), \label{eq:hfgfg1} \\
| f(z)-z| &\leq & \hat{f}(|z|) -|z|.   \label{eq:fgmzfg1}
\end{eqnarray}
\end{lem}
{\sl{Proof}}. We first prove (\ref{eq:hfgfg1}).
There exists complex number $H_n$, which depends on $a_2,\ldots, a_n , b_2,\ldots , b_n$, such that 
\begin{equation*}
f(g(z)) = z+ \sum_{n=2}^{\infty} H_n(a_2,a_3,\ldots , a_n, b_2 , b_3, \ldots, b_n) z^n. 
\end{equation*}
The inequality
\[
|H_n(a_2,a_3,\ldots , a_n, b_2 , b_3, \ldots, b_n)|  \leq  H_n(|a_2|,|a_3|,\ldots , |a_n|, |b_2| , |b_3|, \ldots, |b_n|)
\]
yields
\begin{eqnarray*}
{\widehat{ f \circ g }} (|z|) & \leq & |z|+\sum_{n=2}^{\infty} H_n(|a_2|,|a_3|,\ldots , |a_n|, |b_2| , |b_3|, \ldots, |b_n|) |z|^n \\
                              & = &  \hat { f} (\hat{g}(|z|)).
\end{eqnarray*}
Inequality $(\ref{eq:fgmzfg1})$ is immediate.
\hfill $\Box$ \\
\begin{lem}\label{inequality1}  
Let $d$ be an integer and let $c_{n,r} (n=d,d+1,\ldots ,r =2,3,\ldots)$ be complex numbers such that
\[
f_n(z) := z+ \sum_{r=2}^{\infty} c_{n,r} z^r 
\]
are entire functions for all $n\geq d$.
We set 
\[
F_m(z):= \R_{n=d}^{m} \widehat{f_n} (z) .
\]
Then, for any $z \in \mathbb{C}$ and any integer $m\geq d$ ,
\begin{eqnarray}
\widehat{F_m}(|z|) &\leq & \left( \R_{n=d}^{m} \widehat{f_n} (z) \right) \circ |z|, \label{ineq:fmz1}\\
| \R_{n=d}^{m}f_n(z) -z | &\leq & \left( \R_{n=d}^{m} \widehat{f_n}(z) \right) \circ |z|-|z| , \label{ineq:fmz2} \\
  |\R_{n=d}^{m}f_n(z) | &\leq & \left( \R_{n=d}^{m} \widehat{f_n}(z) \right) \circ |z|. \label{ineq:fmz3}
\end{eqnarray}
\end{lem}
{\sl{Proof}}. We first prove $(\ref{ineq:fmz1})$ by induction on $m$.
If $m=d$, then inequality $(\ref{ineq:fmz1})$ immediately follows. 
Next we suppose that  $(\ref{ineq:fmz1})$ holds for $m (\geq d)$. Then combining $(\ref{eq:hfgfg1})$ in Lemma \ref{ineq:firstlem1} and the inductive assumption yields
\begin{eqnarray*}
\widehat{F_{m+1}}(|z|)  \leq  \widehat{F_{m}}(\widehat{f_{m+1}}(|z|)) \leq  \left( \R_{n=d}^{m} \widehat{f_n} (z) \right) \circ |z| \circ \widehat{f_{m+1}}(|z|) = \left ( \R_{n=d}^{m+1} \widehat{f_n} (z) \right ) \circ |z|.  
\end{eqnarray*}
This completes the proof of  $(\ref{ineq:fmz1})$ .
Next we shall prove $(\ref{ineq:fmz2})$.
\begin{eqnarray*}
|\R_{n=d}^{m}f_n(z) -z| & \leq & \widehat{F_m}( |z|) - |z| \,\,\,\,\, ( \, (\ref{eq:fgmzfg1}) \,\, \text{in Lemma} \ref{ineq:firstlem1}  )\\ 
 & \leq & \left( \R_{n=d}^{m} \widehat{f_n}(z) \right) \circ |z|-|z| \,\, (\, (\ref{ineq:fmz1}) \,\, \text{in Lemma}  \ref{inequality1} ) .
\end{eqnarray*}
This completes the proof of  $(\ref{ineq:fmz2})$. Finally, inequality $(\ref{ineq:fmz3})$ follows from the triangle inequality and $(\ref{ineq:fmz2})$. 
\hfill $\Box$ \\
\begin{lem}\label{inequality2}  
Let $d$ be integer and let $c_{n,r} (n=d,d+1,\ldots ,r =2,3,\ldots)$ be complex numbers such that
\[
f_n(z) := z+ \sum_{r=2}^{\infty} c_{n,r} z^r 
\]
are entire functions.
Suppose that 
\[
C_n := \max_ { r=2,3,\ldots } \{ |c_{n,r}|^{\frac{1}{r-1}}\} > 0 \,\, \text{ for all }\,\, n \geq  d 
\]
Then, for $ |z| < 1/\sum_{n=d}^{m} C_n$ and for any integer $m, d$ with $m\geq d$,
\[
\left( \R_{n=d}^{m} \widehat{f_n}(z) \right) \circ |z| \leq  \frac{|z|}{1-|z|\sum_{n=d}^{m} C_n} .
\]
\end{lem}
{\sl{Proof.}}
We shall prove this by induction on $m$.
Let $m=d$. It follows from $|c_{d,p}| \leq {C_d}^{p-1}$ that, for $ |z| < 1/C_d ,$
\begin{eqnarray*}
\widehat{f_d}(|z|) & = & (z+|c_{d,2}| z^2 + |c_{d,3}| z^3 +\cdots + |c_{d,p}|z^p+ \cdots )\circ |z| \\
   & \leq & |z|+C_d |z|^2+{C_d}^2 |z|^3+ \cdots + {C_d}^{p-1} |z|^p+ \cdots \\
   & = & \frac{|z|}{1-C_d |z|} \, .
\end{eqnarray*}  
Next suppose that the statement of Lemma \ref{inequality2} is valid if $m=N \,(N\geq d)$.\\
Noting
\begin{eqnarray*}
 \frac{1}{{\sum_{n=d}^{N+1}C_n}} \leq  \frac{1}{C_d+C_{N+1}}  < \frac{1}{ C_{N+1}},    
\end{eqnarray*}
we have, for $|z| < ( \sum_{n=d}^{N+1}C_n)^{-1}$,
\begin{eqnarray}
\widehat{f_{N+1}}(|z|) & = & |z|+ |c_{N+1,2}| |z|^2 + |c_{N+1,3}||z|^3 + \cdots + |c_{N+1,p}||z|^p + \cdots \nonumber \\
& \leq & |z| + C_{N+1} |z|^2 + {C_{N+1}}^2 |z|^3 + \cdots + {C_{N+1}}^{p-1} |z|^p + \cdots  \nonumber\\
& = & \frac{|z|}{1-C_{N+1}|z|} . \label{ineq:this112}
\end{eqnarray}
Therefore we have 
\begin{equation}
\widehat{f_{N+1}}(|z|) \leq  \frac{|z|}{1-C_{N+1}|z|} \,\,\, \text{for}\,\,\, |z| < \left( \sum_{n=d}^{N+1}C_n\right)^{-1}.\label{eq:inq1}
\end{equation}
Now we set $ g(z) = z/(1-C_{N+1} z).$
Then $g(z)$ is steadily increasing for $0 \leq z <1/C_{N+1}$.
Noting this and
\begin{eqnarray*}
 \frac{1}{{\sum_{n=d}^{N+1}C_n}}  <  \frac{1}{C_{N+1}},
\end{eqnarray*}
we obtain, for $|z| <1/ \sum_{n=d}^{N+1}C_n$,
\begin{eqnarray}
 \frac{z}{1-C_{N+1}z} \circ |z| & < &  \frac{z}{1-C_{N+1}z} \circ \frac{1}{{\sum_{n=d}^{N+1}C_n}} \nonumber \\
     & = & \frac{1}{{\sum_{n=d}^{N+1}C_n}-C_{N+1}}  =  \frac{1}{{\sum_{n=d}^{N}C_n}} . \label{ineq:this11}
\end{eqnarray}
Combining this with (\ref{eq:inq1}) yields
\[
 \widehat{f_{N+1}}(|z|)  < \frac{1}{{\sum_{n=d}^{N} C_n}} \,\,\, \text{for} \,\,\,|z| < \left( \sum_{n=d}^{N+1}C_n \right)^{-1}.
\]
Hence we can use the inductive assumption as follows:\\
For $|z| <(\sum_{n=d}^{N+1}C_n)^{-1}$,
\begin{eqnarray}
\left( \R_{n=d}^{N} \widehat{f_n}(z) \right) \circ \widehat{f_{N+1}}(|z|) \leq \frac{z}{1-z\sum_{n=d}^{N}C_n} \circ \widehat{f_{N+1}}(|z|) .\label{ineq:im11}
\end{eqnarray}
We set $ h(z) = z/(1-z\sum_{n=d}^{N}C_n).$
Then $h(z)$ is steadily increasing for $ 0 \leq z < 1/\sum_{n=d}^{N}C_n .$
Besides, it follows from (\ref{ineq:this112}) and (\ref{ineq:this11}) that
\begin{eqnarray*}
\widehat{f_{N+1}}(|z|)  \leq  \frac{|z|}{1-C_{N+1}|z|}  <  \frac{1}{{\sum_{n=d}^{N}C_n}} 
\end{eqnarray*}
for $|z|< \left( \sum_{n=d}^{N+1}C_n \right)^{-1}$.
From these, inequality (\ref{ineq:im11}) is rewritten as
\[
\left( \R_{n=d}^{N+1} \widehat{f_n}(z) \right) \circ |z| \leq \frac{z}{1-z\sum_{n=d}^{N}C_n} \circ \frac{|z|}{1-C_{N+1}|z|} 
\]
for $|z|< \left( \sum_{n=d}^{N+1}C_n \right)^{-1}$.
Simplifying the right hand side, we have 

\[
\left( \R_{n=d}^{N+1} \widehat{f_n}(z) \right) \circ |z|\leq \frac{z}{1-z\sum_{n=d}^{N+1}C_n}\circ |z|
\]
for $|z| < 1/\sum_{n=d}^{N+1}C_{n} $.
Hence the induction is complete.
\hfill $\Box$

\begin{lem}\label{cauchy1} 
Let $d$ be integer and let $c_{n,r} (n=d,d+1,\ldots ,r =2,3,\ldots)$ be complex numbers such that
\[
f_n(z) := z+ \sum_{r=2}^{\infty} c_{n,r} z^r 
\]
are entire functions.
Let $C_n$ be the constant given in Lemma \ref{inequality2}, and let $C_n>0$ for every positive integer $n$. Suppose further that $\alpha := \sum_{n=1}^{\infty} C_n$ is convergent.
Moreover we set
\begin{eqnarray*} 
F_N(z):&=& \R_{n=1}^{N} f_n(z).
\end{eqnarray*}
Then, for $|z| \leq 1/(4\alpha)$ and any integers $ N,M$ with $ N > M\geq 1$,
\begin{eqnarray*}
|F_N(z)-F_M(z)| \leq  \frac{1}{\alpha ^2}\sum_{n=M+1}^{N}C_n.
\end{eqnarray*}
\end{lem}
{ \sl{Proof}.}
We set 
\begin{eqnarray*}
 y(z) :&=& \R_{n=M+1}^{N}  f_n(z)  . 
\end{eqnarray*}
Our task is to estimate $ \left| F_N(z)-F_M(z) \right|$.
First note that
\begin{eqnarray*}
\left| F_N(z)-F_M(z) \right| &=& \Bigl| \R_{n=1}^{M} f_{n}(z)   \circ \R_{n=M+1}^{N} f_n(z)   - \R_{n=1}^{M} f_n(z) \Bigl| \\
&=& \left| F_M(y(z))-F_M(z) \right|.
\end{eqnarray*}
To estimate the last expression, we show that
\begin{equation}
 |y(z)| \leq \frac{1}{3\alpha} \,\,\, \text{for} \,\,\, |z|\leq \frac{1}{4\alpha} \label{ineq:yal1}.
\end{equation}
Since 
\begin{equation}
 \frac{1}{4\alpha }= \frac{1}{4\sum_{n=1}^{\infty}C_n} < \frac{1}{\sum_{n=M+1}^{N}C_n}, \label{eq:afi1}
\end{equation}
we can use Lemma \ref{inequality2} with $m=N,\,d=M+1$ as follows:
For $|z| \leq 1/(4\alpha)$,
\begin{eqnarray*}
|y(z)| & \leq & \left( \R_{n=M+1}^{N} \widehat{f_{n}}(z) \right) \circ |z| \quad ( (\ref{ineq:fmz3}) \,\, \text{in Lemma}\,\, \ref{inequality1} ) \\
       & \leq & \frac{|z|}{1-|z|\sum_{n=M+1}^{N}C_n} \quad ( \text{Lemma} \,\, \ref{inequality2} )\\
       & \leq & \frac{1/(4\alpha )}{1-\frac{1}{4\alpha }\sum_{n=M+1}^{N}C_n} \leq \frac{1}{4\alpha -\sum_{n=1}^{\infty}C_n} = \frac{1}{3\alpha }.     
\end{eqnarray*}
Hence we obtain (\ref{ineq:yal1}).
Next, let $\gamma $ be a circle with center $0$ and radius $1/(2\alpha) $.
For $|z| \leq 1/(4\alpha)$, we can use Cauchy's theorem as follows:
\begin{eqnarray}
|F_N(z)-F_M(z)|&=& |F_M(y(z))-F_M(z)| \nonumber\\
                &=& \left| \frac{1}{2\pi i}\int_{\gamma } \left( \frac{F_M(\zeta )}{\zeta -y(z)}-\frac{F_M(\zeta )}{\zeta -z} \right) d\zeta \right|  \,\, (\text{note :}\,\, |y(z)|\leq \frac{1}{3\alpha }) \nonumber \\
                &\leq & \frac{1}{2\pi }\int_{\gamma } \frac{|y(z)-z||F_M(\zeta )|}{|\zeta -y(z)||\zeta-z|} |d\zeta| \nonumber \\
                &\leq &\frac{1}{2\alpha }\frac{|y(z)-z|\max_{|\zeta|=\frac{1}{2\alpha }}|F_M(\zeta )|}{(\frac{1}{2\alpha } - \frac{1}{3\alpha })(\frac{1}{2\alpha }-\frac{1}{4\alpha })} \nonumber \\
                &=& 12 \alpha |y(z)-z| \max_{|\zeta|=\frac{1}{2\alpha }}|F_M(\zeta )| \label{eq:itiyou1}.
\end{eqnarray} 
We first estimate $|y(z) -z|$. For $|z| \leq 1/(4\alpha)$, 
\begin{eqnarray}
|y(z)-z| &\leq & \left( \R_{n=M+1}^{N} \widehat{f_{n}}(z) \right) \circ |z| -|z| \quad ((\ref{ineq:fmz2}) \, \text{in Lemma} \, \ref{inequality1})  \nonumber \\
         &\leq & \frac{|z|}{1-|z|\sum_{n=M+1}^{N}C_n} -|z| \,\,\,\,\,\,\,\,\,(\text{Lemma} \,\, \ref{inequality2} \,\,\  (\text{note}\,\,\, (\ref{eq:afi1})\,)\,) \nonumber \\
         & = & \frac{|z|^2\sum_{n=M+1}^{N}C_n}{1-|z|\sum_{n=M+1}^{N}C_n} \nonumber \\
         &\leq & \frac{1}{16\alpha ^2}\frac{1}{1-\frac{1}{4\alpha }\alpha }\sum_{n=M+1}^{N}C_n = \frac{1}{12\alpha ^2}\sum_{n=M+1}^{N}C_n. \label{eq:itiyou2}      
\end{eqnarray}
Next we shall estimate $\,\max_{|\zeta|=\frac{1}{2\alpha }}|f_M(\zeta )|$.
For $|z| \leq 1/(4\alpha)$, 
\begin{eqnarray}
\max_{|\zeta| =  \frac{1}{2\alpha }}|F_M(\zeta )| &\leq &  \left( \R_{n=1}^{M} \widehat{f_n}(z) \right) \circ \frac{1}{2\alpha } \,\,\, ( (\ref{ineq:fmz3})\,\, \text{in Lemma} \,\,\, \ref{inequality1} ) \nonumber \\
& \leq &  \frac{\frac{1}{2\alpha}}{1-\frac{1}{2\alpha}\alpha} \,\,\,\,\,\,\,\,\, (\text{Lemma} \ref{inequality2})  \nonumber \\
&=& \frac{1}{\alpha} . \label{eq:itiyou3}        
\end{eqnarray}
From (\ref{eq:itiyou1}), (\ref{eq:itiyou2}), and (\ref{eq:itiyou3}), we deduce that
\[
 \left| F_N(z) -F_M(z) \right| \leq \frac{1}{\alpha ^2} \sum_{n=M+1}^{N} C_n  \quad \text{for} \quad |z| \leq \frac{1}{2\alpha }.
\]
This completes the proof of the lemma.
\hfill $\Box$ \\
{\sl{ Proof of Theorem \ref{entire1}.}}
We apply the same notation as in Lemma \ref{cauchy1}. Now we consider two cases.\\
{\bf{Case 1}}  Suppose that there exists a number $m$ such that $C_n =0$ for all $n\geq m$.\\
Then for all $ N > m$,
\begin{eqnarray*}
F_N(z)= \R_{n=1}^{N} f_n(z) = \R_{n=1}^{m-1} f_n(z) \circ \R_{n=m}^{N} f_n(z) = \R_{n=1}^{m-1} f_n(z) .
\end{eqnarray*}
 Accordingly, in this case, the theorem is true.\\
{\bf{Case 2}} Suppose that $C_n>0$ for infinitely many $n$. If there are numbers $n$ such that $C_n=0$, then we can ignore them. Because $C_n=0$ means $f_n(z) = z $, and $z$ is the unit element of composition.
Hence we can suppose without loss of generality that $C_n >0 $ for every positive integer $n$.

From now on, we change the assumption of Theorem \ref{entire1} into
\[
 C_n>0 \,\,\, \text{for all} \,\,\,\, n,
\]
and we shall prove the theorem.

Let $r_1>0$ be any real number. Then it is sufficient to prove that
\[
F_m(z) = \R_{n=1}^{m} f_n(z) 
\]
is uniformly convergent for $|z|\leq r_1$.\\
From Lemma \ref{cauchy1}, we have
\[
F_N(z) = \R_{n=1}^{N} f_n(z) \quad \text{is uniformly convergent on} \,\, |z| \leq \frac{1}{4\alpha}.
\]
In the same way, let $\alpha _m:= \sum_{n=m}^{\infty} C_n $, we have
\[
\R_{n=m}^{N} f_n(z) \quad \text{is uniformly convergent on} \,\, |z| \leq \frac{1}{4\alpha_m }. 
\]
From the assumption, the series $\alpha = \sum_{n=1}^{\infty} C_n$ is convergent,
and hence, for any $ r_1 > 0$, there exists positive integer $m_1$ such that
\begin{equation}
\alpha_{m_1} = \sum_{n=m_1}^{\infty} C_{n} \leq \frac{1}{4r_1} \quad \text{namely} \quad r_1 \leq \frac{1}{4\alpha_{m_1}}  \label{eq:sumfra}. 
\end{equation}
Therefore we find that for any $r_1>0$, there exists number $m_1>1$ such that
\[
\R_{n=m_1}^{N} f_n(z) \quad \text{is uniformly convergent on} \,\, |z| \leq r_1. 
\]
Since the function $\R_{n=1}^{m_1-1}f_n(z)$ is an entire function, we deduce that for any $r_1>0$,
\[
F_N(z)=\R_{n=1}^{m_1-1}f_n(z) \circ \R_{n=m_1}^{N} f_n(z) \quad \text{is uniformly convergent on} \,\, |z| \leq r_1
\]
(see Remark below).
Since $r_1>0$ is an arbitrary real number and $\{F_m(z)\}$ is a sequence of entire functions,
the limit function
\[
\lim_{N \rightarrow \infty} F_N(z) = \R_{n=1}^{\infty} f_n(z)
\]
is also an entire function of $z$.\\
Remark : Let $p(z)$ be an entire function, and let $\{F_m(z)\}_{m=1}^{\infty}$ be a sequence of entire functions. 
Suppose that $F_m(z)$ is uniformly convergent for $|z| \leq r_2$.
Then $p(z) \circ F_m(z) = p(F_m(z))$ is uniformly convergent for $|z| \leq r_2$.
\hfill $\Box$ 
\section{Proof of Proposition \ref{specialvalue}}
In this section, we give a proof of Proposition \ref{specialvalue}. We first prove some lemmas.
\begin{lem}\label{spvlem1}
We set
\[
h_1(z) = \frac{1}{2}(e^{2z}-1),\,\,
h_2(z) = \sin \left( \frac{2z}{\sqrt{3}}+\frac{\pi}{6} \right)-\frac{1}{2},\,\, 
h_3(z) = \left( \sinh \sqrt{z} \right)^2.
\]
Then 
\begin{eqnarray*}
h_1(2z) &=& 2( h_1(z) + h_1(z)^2),\,\,\,\,\,\,\,\, h_1(0)=0,\,\, h_1'(0)=1 \\
h_2(-2z) &=& -2( h_2(z) + h_2(z)^2),\,\,\, h_2(0)=0, \,\, h_2'(0)=1 \\
h_3(4z) &=& 4( h_3(z) + h_3(z)^2),\,\,\,\,\,\,\,\, h_3(0)=0,\,\, h_3'(0)=1 .
\end{eqnarray*}
\end{lem}
{\bf{ \sl{Proof}.}}
The lemma follows from elementary calculations.
\hfill $\Box$

\begin{lem}\label{spvlem2}
Let $s$ be a fixed complex number with $|s|>1$.
Suppose that an entire function $f(z)$ satisfies
\[
f(sz) = s(f(z) +f(z)^2), 
\] 
$f(0)=0$, and $f'(0) =1$. Then
\[
f(z) =  \lim_{N \rightarrow \infty} \R_{n=1}^{N} \left( z+ \frac{z^2}{s^n} \right).
\]
\end{lem}
Hence the function satisfying the above condition is uniquely determined. A generalization of this statement is known (It has been shown in \cite[pp66-68]{Steinmetz} that the inverse function of $ f(z)$ is uniquely determined ).
To prove Lemma \ref{spvlem2}, we require the following lemma.
\begin{lem}\label{dbc}
Suppose that the functions $ p_n(z) , q_n(z), n=1,2,3,\ldots$ are entire functions.
Suppose further that
\[
p_n(z) \,\, \text{is uniformly convergent to} \,\, p(z) \,\, \text{on every compact subset of} \,\,\, \mathbb{C},
\]
and 
\[ 
\quad \lim_{n \rightarrow \infty} q_n(z) = q(z)
\]
exists and becomes entire function.
Then for each $z \in \mathbb{C}$, we have
\[
\lim_{n \rightarrow \infty} p_n(q_n(z)) = p(q(z)) .
\]
\end{lem}
{\sl{Proof.}} 
The lemma follows from
\[
|p(q(z))-p_n(q_n(z))| \leq |p(q(z))-p(q_n(z))| + |p(q_n(z))-p_n(q_n(z))|.
\]
{\bf{ \sl{Proof of Lemma \ref{spvlem2}}.}} 
Let $|s|>1$.
It follows from the assumption that $f(z)$ can be written as in the form
\begin{equation}
f(z) =  z+ \sum_{n=2}^{\infty} a_n(s) z^n \label{fsif3},
\end{equation}
and satisfies
\[
f(z) = s \left( z+ z^2 \right) \circ f\left( \frac{z}{s} \right) .
\]
We repeatedly use the relation as follows:
\begin{eqnarray*}
f(z) &=& \left( z+ \frac{z^2}{s} \right) \circ \left( sf\left( \frac{z}{s} \right) \right)  \\
     &=& \left( z+ \frac{z^2}{s} \right) \circ sz \circ  s(z+z^2)\circ f\left(\frac{z}{s^2}\right) \\
     &=& \left( z+ \frac{z^2}{s} \right) \circ s^2z \circ  (z+z^2) \circ \frac{z}{s^2} \circ s^2z \circ f\left(\frac{z}{s^2}\right) \\
     &=& \left( z+ \frac{z^2}{s} \right) \circ  \left( z+\frac{z^2}{s^2} \right) \circ \left( s^2f\left(\frac{z}{s^2}\right) \right) \\
     &=& \cdots \\
     &=& \left( \R_{n=1}^{N} \left( z+ \frac{z^2}{s^n} \right) \right) \circ  \left( s^Nf\left(\frac{z}{s^N}\right) \right). 
\end{eqnarray*}
Hence we deduce 
\begin{equation}
f(z) = \lim_{N \rightarrow \infty} \left\{ \left( \R_{n=1}^{N} \left( z+ \frac{z^2}{s^n} \right) \right) \circ  \left( s^Nf\left(\frac{z}{s^N}\right) \right) \right\} . \label{eq:rlr1}
\end{equation}
From $(\ref{fsif3})$, we have
\begin{equation}
 \lim_{N \rightarrow \infty} s^N f\left(\frac{z}{s^N}\right) = z \quad \text{for every} \,\, z \in \mathbb{C}. \label{eq:linz}
\end{equation}
Since $\R_{n=1}^{N} \left( z+ \frac{z^2}{s^n} \right)$ is convergent on every compact subset of $\mathbb{C}$ from Theorem \ref{entire1}, combining Equation (\ref{eq:linz}), Equation (\ref{eq:rlr1}), and Lemma \ref{dbc} yields
\[
f(z) =  \lim_{N \rightarrow \infty} \R_{n=1}^{N} \left( z+ \frac{z^2}{s^n} \right).
\]
\hfill $\Box$\\
{\sl{Proof of Proposition \ref{specialvalue}.}}
The proposition follows by Lemma \ref{spvlem1} and Lemma \ref{spvlem2}.
\hfill $\Box$

{{\sl{Acknowledgments}}}:
The author is very grateful to Professor Akio Fujii, Professor Noboru Aoki, Professor Shigeki Akiyama,  Assistant Professor Tomoki Kawahira, and Assistant Professor Masatoshi Suzuki for their useful advice and their encouragement. In particular, Assistant Professor Tomoki Kawahira taught the author Lemma \ref{ineq:firstlem1} and Lemma \ref{inequality1}, and Professor Akio Fujii and Assistant Professor Masatoshi Suzuki offered  many helpful comments.

\end{document}